**Conditional Maximum Lq-Likelihood Estimation for Regression Model with Autoregressive Error Terms**

Y. Güney, Y. Tuaç, Ş. Özdemir and O. Arslan

Ankara University, Faculty of Science, Department of Statistics, 06100 Ankara/Turkey

Afyon Kocatepe University, Faculty of Science, Department of Statistics, Afyonkarahisar, Turkey

ydone@ankara.edu.tr, ytuac@ankara.edu.tr, senayozdemir@aku.edu.tr, oarslan@ankara.edu.tr

**Abstract**

In this article, we consider the parameter estimation of regression model with $p^{th}$ order autoregressive (AR(p)) error term. We use the Maximum Lq-likelihood (MLq) estimation method that is proposed by Ferrari and Yang (2010a), as a robust alternative to the classical maximum likelihood (ML) estimation method to handle the outliers in the data. After exploring the MLq estimators for the parameters of interest, we provide some asymptotic properties of the resulting MLq estimators. We give a simulation study and a real data example to illustrate the performance of the new estimators over the ML estimators and observe that the MLq estimators have superiority over the ML estimators when outliers are present in the data.

**Keywords:** autoregressive stationary process; conditional maximum Lq-likelihood; linear regression.

1. **Introduction**

Incorporating the autoregressive error terms into the linear regression models is a useful way for analyzing the relationships between economic indicators and has also been attracted a good deal interest in both statistical theory and applications. In literature, classical parameter estimation methods are generally used to estimate the parameters of the regression model with autoregressive error terms (AR(p)). For instance, Cochrane and Orcutt (1949) have considered some modification of the ordinary least squares (OLS) estimation method for autoregressive error terms regression model. Beach and Mackinnon (1978) have used maximum likelihood (ML) estimation method to estimate the parameters of AR(1) error term regression model. Alpuim and El-Shaarawi (2008) have estimated the parameters of the regression model with AR(p) error term using the OLS estimation method. They have also used the ML estimation and conditional maximum likelihood (CML) estimation method under the assumption of normality and studied the asymptotic properties of the resulting estimators. Tuac et al. (2017) have considered linear regression model with AR(p) errors with Student's *t*-distribution as a heavy tailed alternative to the normal distribution and used CML estimation method to obtain the model parameters.

Under normality assumption, CML and OLS estimation methods are commonly used estimation procedures for the regression model with autoregressive error terms. Although they are appropriate choices in estimating the parameters of regression with autoregressive errors, they are highly sensitive to outliers. To handle this problem, we will use the MLq estimation method proposed by Ferrari and Yang (2010a) to estimate the parameters of the AR(p) error term regression model.

MLq estimation method have recently proposed and have gained considerable attention in the past decade. For example, Ferrari and Paterlini (2009) have investigated the behavior of the MLq estimation on both simulated data and on real-world time series for extreme quantile estimation. Ferrari and Paterlini (2010b) have applied the MLq estimation to expected return and volatility estimation of financial asset returns under multivariate normality. Huang, Lin and Ren (2013) have proposed a generalized form of the classical likelihood ratio statistic by using the Lq -likelihood ratio (LqR) statistic based on the MLq estimation for hypothesis testing problem for the shape parameter of the GEV distribution and showed that the asymptotic behavior of proposed statistic characterize by the degree of tuning parameter q. Qin and Priebe (2013) have proposed a new EM algorithm namely an expectation-maximization algorithm with Lq-likelihood (EMLq) which addresses MLq estimation within the EM framework for mixture models. Qin and Priepe (2016) have introduced a robust hypothesis testing procedure: the Lq-likelihood-ratio-type test (LqRT). By deriving the asymptotic distribution of this test statistic, the authors have demonstrated its robustness both analytically and numerically, and they investigated the properties of both its influence function and its breakdown point. Also Ozdemir et al. (2019) use the MLq estimation method to estimate the parameters of Marshall-Olkin extended Burr XII distribution and show that MLq estimation method outperform the ML. Recently, Dogru et al. (2018) propose parameter estimation of the multivariate t distribution using the MLq estimation, provide that unlike the ML estimation the degrees of freedom parameters can be estimated along with the other parameters, and still gain the robustness.

In this paper, we consider the conditional maximum Lq-likelihood (CMLq) estimation method for the autoregressive error terms regression models under normality assumption. We obtain the parameter estimation for all the parameters. We give an extensive simulation study to compare the performances of the CML and the CMLq estimation methods. The performance of the two sets of estimators is evaluated for different data structure including outlier cases. The simulation results show that the CMLq estimators can reduce the effects of the outliers if the tuning parameter $q$ is less than one. Note that q is a key parameter to maintain the robustness in MLq estimation method. It is considered as a robustness tuning parameter and choosing it very important. In this paper we choose q that minimizes robust AIC. Further we provide the asymptotic distribution of the proposed estimator and use the asymptotic covariance matrix to form the asymptotic confidence intervals for the parameters.

The rest of the paper is organized as follows. In the next section, we briefly describe the regression models with AR(p) error terms and the CML estimation method for the parameters of interest. A brief description of the MLq estimation method is summarized in Section 3. Then we obtain the CMLq estimators for the parameters of the regression models with AR(p) error terms in Section 4. Since the estimators cannot be obtained in explicit form the iteratively reweighted algorithm steps given in Section 5. Section 6 presents the asymptotic normality of the proposed estimators. In Section 7, Monte Carlo simulation study and a real data example are presented to compare the performance of CMLq and the CML estimation methods in terms of RMSE at different data structure scenarios including outliers. Finally, some concluding remarks are given in Section 8.

## 2. Regression Models with Autoregressive Error

Consider the linear regression model with the error terms follow a stationary AR(p) process

$$y_t = \sum_{i=1}^{M} x_{t,i} \beta_i + e_t , \tag{1}$$

$$e_t = \phi_1 e_{t-1} + \cdots + \phi_p e_{t-p} + a_t, \quad t = 1,2,\ldots,N \tag{2}$$

where, $y_t$ is the response variable, $x_{t,i}$ are explanatory variables, $\beta_i$ are unknown regression parameters and $\phi_j$ are unknown autoregressive parameters. The $a_t$ is normally and independently distributed such that $E(a_t) = 0$ and $Var(a_t) = \sigma^2$.

To simplify the notation, we denote $a_t = \Phi(B)e_t$ and here, $B$ is called the Backshift operator. By doing so, an alternative form of the model (1) is

$$\Phi(B)y_t = \sum_{i=1}^{M} \beta_i \Phi(B)x_{t,i} + a_t, \quad t = p+1, p+2, \ldots, N, \tag{3}$$

where

$$\Phi(B)y_t = y_t - \phi_1 y_{t-1} - \cdots - \phi_p y_{t-p}, \tag{4}$$

$$\Phi(B)x_{t,i} = x_{t,i} - \phi_1 x_{t-1,i} - \cdots - \phi_p x_{t-p,i} \tag{5}$$

$\phi_j$ are the AR(p) model parameters for $j = 1,2,\ldots,p$.

Since the exact likelihood function could be well approximated by the conditional likelihood function (Ansley, 1979) we will first give the CML estimation method which are used mainly in cases where ML estimates are difficult to compute.

## 2.1 Parameter Estimation of AR(p) Error Terms Regression Model with CML

Under the assumptions associated with equation (3) the conditional log-likelihood function will be as follows (Alpuim and El-Shaarawi, 2008).

$$lnL = c - \frac{N}{2}ln\sigma^2 - \frac{1}{2\sigma^2}\sum_{t=p+1}^{N}\left(\Phi(B)y_t - \sum_{i=1}^{M}\beta_i\,\Phi(B)x_{t,i}\right)^2 \tag{6}$$

To obtain the estimating equations, the derivatives of the conditional log-likelihood function with respect to unknown parameters are taken and set to zero. The derivatives of the conditional log-likelihood function are given in Appendix A.

Rearranging the estimating equations, we get the following forms of estimators.

$$\underline{\hat{\beta}} = \left[\sum_{t=p+1}^{N}\widehat{\Phi}(B)x_t\widehat{\Phi}(B)x_t^T\right]^{-1}\left[\sum_{t=p+1}^{N}\widehat{\Phi}(B)y_t\widehat{\Phi}(B)x_t\right] \tag{7}$$

$$\underline{\hat{\phi}} = \boldsymbol{R}^{-1}\left(\underline{\hat{\beta}}\right)R_0\left(\underline{\hat{\beta}}\right) \tag{8}$$

$$\hat{\sigma}^2 = \frac{1}{N-p}\sum_{t=p+1}^{N}\left(\widehat{\Phi}(B)y_t - \underline{\hat{\beta}}\widehat{\Phi}(B)x_t\right)^2 \tag{9}$$

where

$$R_0\left(\underline{\beta}\right) = \begin{bmatrix} \sum_{t=p+1}^{N}e_t e_{t-1} \\ \sum_{t=p+1}^{N}e_t e_{t-2} \\ \vdots \\ \sum_{t=p+1}^{N}e_t e_{t-p} \end{bmatrix}, \boldsymbol{R}\left(\underline{\beta}\right) = \begin{bmatrix} \sum_{t=p+1}^{N}e_{t-1}^2 & \sum_{t=p+1}^{N}e_{t-1}e_{t-2} & \cdots & \sum_{t=p+1}^{N}e_{t-1}e_{t-p} \\ \sum_{t=p+1}^{N}e_{t-2}e_{t-1} & \sum_{t=p+1}^{N}e_{t-2}^2 & \cdots & \sum_{t=p+1}^{N}e_{t-2}e_{t-p} \\ \vdots & \vdots & \ddots & \vdots \\ \sum_{t=p+1}^{N}e_{t-p}e_{t-1} & \sum_{t=p+1}^{N}e_{t-p}e_{t-2} & \cdots & \sum_{t=p+1}^{N}e_{t-p}^2 \end{bmatrix} \tag{10}$$

and $\widehat{\Phi}(B)$ is the backshift operator with the estimates of $\phi_j$.

Equations (7) and (9) can be written in vector form as follows

$$\underline{\hat{\beta}} = \left[\widehat{\Phi}(B)X^T\widehat{\Phi}(B)X\right]^{-1}\left[\widehat{\Phi}(B)X^T\widehat{\Phi}(B)\underline{Y}\right] \tag{11}$$

$$\hat{\sigma}^2 = \frac{1}{N-p}[\widehat{\Phi}(B)\underline{Y} - \widehat{\Phi}(B)X\underline{\hat{\beta}}]^T[\widehat{\Phi}(B)\underline{Y} - \widehat{\Phi}(B)X\underline{\hat{\beta}}] \tag{12}$$

where

$$\widehat{\Phi}(B)X = [\widehat{\Phi}(B)x_{t,i}],$$

$$\widehat{\Phi}(B)Y = [\widehat{\Phi}(B)y_t].$$

These vector forms are important to implement the IRA algorithm to compute the estimates and will be our updating equations.

## 3. Maximum Lq Likelihood Estimation Method

Ferrari and Yang (2010a) introduced the MLq estimation method based on q entropy and defined as follows. Suppose $x = (x_1, x_2, \ldots, x_n)$ be a random sample from a distribution with probability density function $f(x; \theta)$ with $\theta \in \Theta$. The MLq estimator of $\theta$ is defined as

$$\tilde{\theta}_{ML_qE} = \underset{\theta \in \Theta}{\mathrm{argmax}} \sum_{i=1}^{n} L_q(f(x_i; \theta)) \ , \ q > 0 \tag{13}$$

where $L_q : (0, \infty) \to \mathbb{R}$ called non-extensive entropy or *q*-order entropy is defined by

$$L_q(u) = \begin{cases} \log u & , \quad q = 1 \\ u^{1-q} - 1/{1-q} & , \quad otherwise \end{cases}$$

(Havrda and Charvát, 1967 and Tsallis, 1988).

Define

$$U^*(x;\theta,q) = \nabla_\theta L_q(f(x;\theta)) = U(x;\theta)f(x;\theta)^{1-q}.$$

where $U(x;\theta) = \nabla_\theta \log f(x;\theta)$. $\nabla$ denotes derivative operator. Then Lq-likelihood equations have the form

$$\sum_{i=1}^{n} U^*(x_i;\theta,q) = \sum_{i=1}^{n} U(x_i;\theta)f(x_i;\theta)^{1-q} = 0. \tag{14}$$

The MLq estimation method can be regarded as a generalization of the ML estimation method. In particular, when $q$ is equal to 1, MLq and ML estimation methods are equal. It is also easy to see that when $\to 1$, $L_q(u) \to \log u$. In this sense the MLq approaches the ML. Specially, MLq estimator belongs to the class of M-estimators (Hampel, et al. 1986, Huber et al. 2009, Maronna et. al. 2006), as q is fixed (Ferrari and Paterlini, 2009). MLq estimator also minimizes power divergences. The family of power divergences has various special cases for specific values of q (For more details, see Ferrari and Paterlini, 2009). For instance, for q = 1/2 the MLq estimator is minimum Hellinger distance estimator.

Classical weighted likelihood approach has been developed to deal with the disadvantages of maximum likelihood approach by modifying the role of observations by means of the weights. The MLq likelihood score functions are the same as the weighted likelihood score functions. That is why the MLq estimation can be seen as weighted likelihood estimators. However, it is important to note that the MLq estimator isn't obtained with the idea of weighting the likelihood function. The difference from most weighted likelihood estimators is that the weights are a function of pdf. The weight function is also based on tuning parameter $q$. Therefore, choices $q$ is the key for this approach. MLq estimation method give low or high weights depending on $q < 1$ or $q > 1$ to the extreme observations, by doing so it reduces the effects of outliers. In particular, all the observations have the same weight when $q=1$. Taking the $q$ value slightly different makes it possible to obtain better estimates in terms of balance between bias and variance when the sample size is small. MLq estimation method balances two apparently contrasting needs: efficiency and robustness by a proper choice of $q$. It provides strong robustness at expense of a slightly reduced efficiency in presence of outliers. Specially, the effect of extreme observations is reduced by taking $q<1$. On the contrary, when $q>1$, the effect of the observations corresponding to density values close to zero is accentuated (Ferrari and Yang, 2007). A detailed discussion of the role of the tuning parameter $q$ see Ferrari and Yang (2007).

**4. Parameter Estimation of AR(p) Error Terms Regression Model with CMLq**

Let's $(x_i, y_i)$, $i = 1, 2, \dots, n$, be a random sample from model given in (3) with the assumption that $a_t$ is normal and $\underline{\theta}' = (\beta_1, \beta_2, \dots, \beta_M, \phi_1, \phi_2, \dots, \phi_p, \sigma^2)$ be the parameter vector of the model. Then the CMLq estimator of $\underline{\theta}$ is defined as

$$\underline{\tilde{\theta}}_{ML_qE} = \underset{\underline{\theta} \in \Theta}{\operatorname{argmax}} \sum_{t=p+1}^{N} L_q(f(a_t|a_1, a_2, \dots, a_p, \underline{\theta})). \tag{15}$$

Define

$$U(a_t; \underline{\theta}) = \nabla_\theta \log f(a_t; \underline{\theta}),$$

$$U^*(a_t; \underline{\theta}, q) = U(a_t; \underline{\theta}) f(a_t; \underline{\theta})^{1-q}, \quad t = p+1, \dots, N.$$

Then Lq-likelihood equations are

$$\sum_{t=p+1}^{N} U^*(a_t; \underline{\theta}, q) = \sum_{t=p+1}^{N} U(a_t; \underline{\theta}) f(a_t; \underline{\theta})^{1-q} = \underline{0}. \tag{16}$$

The elements of $U(a_t; \underline{\theta})$ and $U^*(a_t; \underline{\theta}, q)$ are given in Appendix A and Appendix B, respectively.

The CMLq estimators of the parameters are the solutions of (16). From these equations we get the following forms of estimators. Provided that $\left[\sum_{t=p+1}^{N} \omega_t(v) \widehat{\Phi}(B) x_t \widehat{\Phi}(B) x_t^T\right]^{-1}$ and $\boldsymbol{R}_\omega^{-1}(\hat{\beta})$ exist.

$$\underline{\tilde{\beta}} = \left[\sum_{t=p+1}^{N} \omega_t \tilde{\Phi}(B)x_t \tilde{\Phi}(B)x_t^T\right]^{-1} \left[\sum_{t=p+1}^{N} \omega_t \tilde{\Phi}(B)y_t \tilde{\Phi}(B)x_t\right], \tag{17}$$

$$\underline{\tilde{\phi}} = \boldsymbol{R_\omega}^{-1}\left(\underline{\tilde{\beta}}\right) R_{\omega 0}\left(\underline{\tilde{\beta}}\right), \tag{18}$$

$$\tilde{\sigma}^2 = \frac{1}{\sum_{t=p+1}^{N} \omega_t} \sum_{t=p+1}^{N} \omega_t \left(\tilde{\Phi}(B)y_t - \tilde{\Phi}(B)x_t^T \underline{\tilde{\beta}}\right)^2 \tag{19}$$

where

$$R_{\omega 0}\left(\underline{\beta}\right) = \begin{bmatrix} \sum_{t=p+1}^{N} \omega_t e_t e_{t-1} \\ \sum_{t=p+1}^{N} \omega_t e_t e_{t-2} \\ \vdots \\ \sum_{t=p+1}^{N} \omega_t e_t e_{t-p} \end{bmatrix},$$

$$\boldsymbol{R_\omega}\left(\underline{\beta}\right) = \begin{bmatrix} \sum_{t=p+1}^{N} \omega_t e_{t-1}^2 & \sum_{t=p+1}^{N} \omega_t e_{t-1} e_{t-2} & \cdots & \sum_{t=p+1}^{N} \omega_t e_{t-1} e_{t-p} \\ \sum_{t=p+1}^{N} \omega_t e_{t-2} e_{t-1} & \sum_{t=p+1}^{N} \omega_t e_{t-2}^2 & \cdots & \sum_{t=p+1}^{N} \omega_t e_{t-2} e_{t-p} \\ \vdots & \vdots & \ddots & \vdots \\ \sum_{t=p+1}^{N} \omega_t e_{t-p} e_{t-1} & \sum_{t=p+1}^{N} \omega_t e_{t-p} e_{t-2} & \cdots & \sum_{t=p+1}^{N} \omega_t e_{t-p}^2 \end{bmatrix}$$

We can obtain the following vector forms related the estimators $\underline{\hat{\beta}}$ and $\hat{\sigma}^2$

$$\underline{\tilde{\beta}} = \left[\tilde{\boldsymbol{\Phi}}(B)\boldsymbol{X}^T \boldsymbol{W} \tilde{\boldsymbol{\Phi}}(B)\boldsymbol{X}\right]^{-1}\left[\tilde{\boldsymbol{\Phi}}(B)\boldsymbol{X}^T \boldsymbol{W} \tilde{\boldsymbol{\Phi}}(B)\underline{Y}\right], \tag{20}$$

$$\tilde{\sigma}^2 = \frac{1}{trace(\boldsymbol{W})}\left[\tilde{\boldsymbol{\Phi}}(B)\underline{Y} - \tilde{\boldsymbol{\Phi}}(B)\boldsymbol{X}\underline{\tilde{\beta}}\right]^T \boldsymbol{W}\left[\tilde{\boldsymbol{\Phi}}(B)\underline{Y} - \tilde{\boldsymbol{\Phi}}(B)\boldsymbol{X}\underline{\tilde{\beta}}\right], \tag{21}$$

where

$$\widetilde{\Phi}(B)X = \left[\widetilde{\Phi}(B)x_{t,i}\right]_{\substack{t=p+1,\dots,N,\\ i=1,\dots,M}}$$

$$\widetilde{\Phi}(B)Y = \left[\widetilde{\Phi}(B)y_t\right]_{t=p+1,\dots,N},$$

$$W = diag\{\omega_t\}_{t=p+1,\dots,N}.$$

From equations (17)-(19), CMLq estimators can be viewed as a weighted version of the CML estimators given in (7)-(9). Here the weights are proportional to the $(1-q)^{th}$ power of the normal distribution density. If $\omega_t = 1$ then (17)-(19) gives the CML estimates of the parameters. Since the weight function $\omega_t$ is a decreasing function of $\left(\Phi(B)y_t - \sum_{i=1}^{M}\beta_i \Phi(B)x_{t,i}\right)^2/\sigma^2$ as $q<1$, the observations with larger residuals receive small weights. When $q$ approaches to zero, the weights get smaller. Thus, the weight function down-weights the effect of the outliers on the estimation procedure. The tuning parameter $q$ balances the efficiency and the robustness of the estimator. When q gets closer to one the resulting estimators get closer to the CML estimators. On the other hand, smaller values of q produce estimators that less sensitive to the outliers but not as efficient as CML.

Concerning the explicit forms of the estimators, we observe that the estimators given in (17) - (19) depend on the weights and the weights are also function of the estimators. Therefore, explicit solutions cannot be obtained from this system of equations. It is also clear from equations (17) - (19) that for the fixed value of $q$, the estimation problem can be solved in terms of a weighting process. In this study, we use the iteratively re-weighted algorithm (IRA) to solve this problem.

## 5. Iteratively Reweighted Algorithm to Compute the Estimates

Let $m \in \{0,1,2,...\}$ denotes the iteration step.

(i) Set the initial values $\underline{\beta}^{(0)}, \underline{\phi}^{(0)}$ and $\sigma^{2(0)}$ and fix a stopping rule ($\varepsilon$).

(ii) Calculate the following weight function for $m = 0,1,2 ...$

$$\omega_t^{(m)} = \left( \frac{1}{\sqrt{2\pi\sigma^{2(m)}}} \exp\left\{ -\frac{\left(\Phi^{(m)}(B)y_t - \sum_{i=1}^{M}\beta_i^{(m)} \Phi^{(m)}(B)x_{t,i}\right)^2}{2\sigma^{2(m)}} \right\} \right)^{1-q}$$

(iii) Calculate $R_\omega\left(\underline{\beta}^{(m)}\right)$ and $R_{\omega 0}\left(\underline{\beta}^{(m)}\right)$

$$\underline{\phi}^{(m+1)} = R_\omega^{-1}\left(\underline{\beta}^{(m)}\right) R_{\omega 0}\left(\underline{\beta}^{(m)}\right)$$

(iv) Using $\omega_t^{(m)}$ and $\underline{\phi}^{(m+1)}$ calculate

$$\underline{\beta}^{(m+1)} = \left[\Phi^{(m+1)}(B)X^T W^{(m)} \Phi^{(m+1)}(B)X\right]^{-1} \left[\Phi^{(m+1)}(B)X^T W^{(m)} \Phi^{(m+1)}(B)\underline{Y}\right]$$

where $W^{(m)} = \text{diag}\{\omega_t^{(m)}\}_{t=p+1}^{n}$.

(v) Using $W^{(m)}, \underline{\phi}^{(m+1)}$ and $\underline{\beta}^{(m+1)}$ calculate

$$(\sigma^2)^{(m+1)} = q \frac{1}{\text{trace}(W^{(m)})} \left[\Phi^{(m+1)}(B)\underline{Y} - \Phi^{(m+1)}(B)X\underline{\beta}^{(m+1)}\right]^T W^{(m)} \left[\Phi^{(m+1)}(B)\underline{Y} - \Phi^{(m+1)}(B)X\underline{\beta}^{(m+1)}\right]$$

(vi) If $\left\|\underline{\beta}^{(m+1)} - \underline{\beta}^{(m)}\right\| < \varepsilon$, $\left\|\underline{\phi}^{(m+1)} - \underline{\phi}^{(m)}\right\| < \varepsilon$ and $\left|(\sigma^2)^{(m+1)} - (\sigma^2)^{(m)}\right| < \varepsilon$ then stop, else repeat the steps (ii-vi) until the convergence condition is satisfied.

## 6. Asymptotic Distribution of CMLq Estimator and Asymptotic Confidence Interval

In this section, the asymptotic covariance matrix of the CMLq estimator of the parameters of the autoregressive error terms regression model under the assumptions given in Section 2 is obtained to construct the asymptotic confidence intervals for the parameters of interest.

Ferrari and Yang (2010a), provide the asymptotic distribution of MLq estimator for the parameters of any exponential family of distributions under some assumptions. We will briefly describe their results. Let $\theta_n^*$ be the value such that

$$E_{\theta_0}\bigl(U^*(\theta_n^*; X, q_n)\bigr) = 0.$$

They state that $\theta_n^* = \frac{\theta}{q_n}$ where $\theta$ is the true parameter. They call $\theta_n^*$ the surrogate parameter of $\theta$. As they point out, since the actual target of $\tilde{\theta}_n$ given in (13) is $\theta_n^*$, $q_n$ must converge to one to obtain the asymptotic unbiasedness of $\tilde{\theta}_n$. Therefore, under the following conditions they show that Lq-likelihood equation has a solution, it is unique and maximizes the Lq-likelihood function in the parameter space.

A1. $q_n > 0$ is a monotone sequence such that $q_n \to 1$ as $n \to \infty$.

A2. The parameter space $\Theta$ is compact and the parameter $\theta$ is an interior point in $\Theta$.

Further, with these assumptions the asymptotic distribution of MLq estimator is

$$\sqrt{n}\mathbf{V}_q^{-1/2}(\tilde{\theta}_n - \theta_n^*) \to N_p(\underline{0}, \mathbf{I}_p),$$

where $\mathbf{I}_p$ is the $(p \times p)$ identity matrix and

$$\mathbf{V}_q(\theta_n^*) = \mathbf{J}(\theta_n^*)^{-1}\mathbf{K}(\theta_n^*)\mathbf{J}(\theta_n^*)^{-1}$$

$$\mathbf{J}(\theta_n^*) = E\left(\nabla_{\underline{\theta}} U^*(\theta_n^*; X, q_n)\right)$$

$$\boldsymbol{K}(\theta_n^*) = E\left(U^*(\theta_n^*; X, q_n)'U^*(\theta_n^*; X, q_n)\right).$$

Here $U^*(\theta_n^*; X, q_n) = U(\theta_n^*; X)f(X, \theta_n^*)^{1-q_n}$ is the first and $\nabla_{\underline{\theta}} U^*(\theta_n^*; X, q_n) = \nabla_{\underline{\theta}}\left(U(\theta_n^*; X)f(X, \theta_n^*)^{1-q_n}\right)$ is the second partial derivatives of Lq-likelihood function.

Concerning the autoregressive error terms regression model under the assumptions given in Section 2 let $\underline{\theta}^{(0)} = \left(\beta_1^{(0)}, \beta_2^{(0)}, \dots, \beta_M^{(0)}, \phi_1^{(0)}, \phi_2^{(0)}, \dots, \phi_p^{(0)}, \sigma^{(0)2}\right)$ be the real parameter vector. Further, let $U(\underline{\theta}; a_t)$ be the score vector. It is known that for all $\underline{\theta} = (\beta_1, \beta_2, \dots, \beta_M, \phi_1, \phi_2, \dots, \phi_p, \sigma^2) \in \Theta$, $E_{\underline{\theta}}[U(\underline{\theta}; a_t)] = \underline{0}$, namely,

$$\int U(\underline{\theta}; a_t)f(a_t, \underline{\theta})d(a_t) = 0. \tag{22}$$

The modified score vector is such that $U^*(\underline{\theta}, a_t, q_n) = U(\underline{\theta}; a_t)f(a_t, \underline{\theta})^{1-q_n}$ and so

$$E_{\underline{\theta}^{(0)}}[U^*(\underline{\theta}; a_t, q_n)] = \int U(\underline{\theta}; a_t)f(a_t, \underline{\theta})^{1-q_n} f(a_t, \underline{\theta}^{(0)})d(a_t).$$

By rearranging the function to be integrate, we obtain

$$E_{\underline{\theta}^{(0)}}[U^*(\underline{\theta}; a_t, q_n)] = \frac{(2\pi)^{\frac{q_n-1}{2}}\sigma^q}{\sigma^{(0)}} \int U(\underline{\theta}; a_t)f(a_t, \underline{\theta}) \frac{exp\left(\frac{-1}{2\sigma^{(0)2}}a_t^{(0)2}\right)}{exp\left(\frac{-q_n}{2\sigma^2}a_t^2\right)} d(a_t) \tag{23}$$

where $a_t^{(0)} = a_t|_{\underline{\theta}=\underline{\theta}^{(0)}}$.

Define $\underline{\theta}_n^* = \left(\beta_1^{(0)}, \beta_2^{(0)}, \ldots, \beta_M^{(0)}, \phi_1^{(0)}, \phi_2^{(0)}, \ldots, \phi_p^{(0)}, q_n \sigma^{(0)2}\right)$ and substitute $\underline{\theta}_n^*$ with $\underline{\theta}$ in equation (23), then we obtain

$$E_{\underline{\theta}^{(0)}}[U^*(\underline{\theta}_n^*; a_t, q_n)] = \frac{(2\pi)^{\frac{q_n-1}{2}}(q_n\sigma^{(0)})^{q_n}}{\sigma^{(0)}} \int U(\underline{\theta}_n^*; a_t) f(a_t, \underline{\theta}_n^*) \frac{\exp\left(\frac{-1}{2\sigma^{(0)2}} a_t^{(0)2}\right)}{\exp\left(\frac{-q_n}{2(q_n\sigma^{(0)2})} (a_t|_{\underline{\theta}=\underline{\theta}^*})^2\right)} d(a_t). \tag{24}$$

Since the parameters which are related to the location in $\underline{\theta}^{(0)}$ and $\underline{\theta}_n^*$ are equal to each other, $a_t|_{\underline{\theta}=\underline{\theta}^*} = a_t|_{\underline{\theta}=\underline{\theta}^{(0)}} = a_t^{(0)}$. Therefore, following equation is provided

$$\frac{\exp\left(\frac{-1}{2\sigma^{(0)2}} a_t^{(0)2}\right)}{\exp\left(\frac{-q_n}{2(q_n\sigma^{(0)2})} (a_t|_{\underline{\theta}=\underline{\theta}^*})^2\right)} = \frac{\exp\left(\frac{-1}{2\sigma^{(0)2}} a_t^{(0)2}\right)}{\exp\left(\frac{-q_n}{2(q_n\sigma^{(0)2})} a_t^{(0)2}\right)} = 1.$$

Then the equation (24) can rewritten as follows.

$$E_{\underline{\theta}^{(0)}}[U^*(\underline{\theta}_n^*; a_t, q_n)] = \frac{(2\pi)^{\frac{q_n-1}{2}}(q_n\sigma^{(0)})^{q_n}}{\sigma^{(0)}} \int U(\underline{\theta}_n^*; a_t) f(a_t, \underline{\theta}_n^*) d(a_t)$$

Knowing that the equation (22) is valid for all $\underline{\theta} \in \Theta$, we obtain

$$\int U(\underline{\theta}_n^*; a_t) f(a_t, \underline{\theta}_n^*) d(a_t) = 0$$

and

$$E_{\underline{\theta}^{(0)}}[U^*(\underline{\theta}_n^*; a_t, q_n)] = 0.$$

Finally, for our problem the surrogate parameter will be $\underline{\theta}_n^* = \left(\beta_1^{(0)}, \beta_2^{(0)}, \ldots, \beta_M^{(0)}, \phi_1^{(0)}, \phi_2^{(0)}, \ldots, \phi_p^{(0)}, q_n \sigma^{(0)2}\right)$. Notice that only the scale parameter $\sigma$ is different from the true parameter value $\sigma^{(0)}$, it is depend on the tuning parameter $q_n$. For this reason, CMLq estimator is expected to make substantial improvements in the prediction of $\sigma$ without much improvement on location parameters.

Also note that, Cavalieri (2002) has shown the same procedure for measurement error models to determine the surrogate parameter.

***Remark.*** *The considered model contains m+p+1 parameter which are $\beta_i, i = 1,2, ..., M$ and $\phi_j$, $j = 1,2, ..., p$ are the model parameter which are related to the location parameter of the distribution. On the other hand $\sigma$ is the scale parameter and unlike the location parameter, it controls the shape of the distribution (the kurtosis of the distribution). Since the MLq method is based on the $q_n^{th}$ power of the pdf, the tuning parameter $q_n$ only affects the shape of the distribution not the location. Therefore it is observed from the asymptotic properties that only scale parameter $\sigma$ is affected from q. For this reason in simulation study, we have redefined $\sigma$ by multiplying the resulting estimates with q.*

By maximizing the Lq-likelihood function we actually target to get $\underline{\theta}_n^*$. Therefore, to get $\underline{\theta}^{(0)}$, $q_n$ should tends to one as n tends to infinity. Therefore, similar to the Ferrari and Yang (2010a) the assumptions A1 and A2 should be hold to get asymptotic unbiasedness of $\underline{\tilde{\theta}}_n$. These assumptions also guarantee that there is a solution of the Lq-likelihood equation, it is a maximizer of the Lq-likelihood function and it is unique.

Similar to Ferrari and Yang (2010a) the asymptotic distribution of the MLq estimator for the autoregressive error terms regression model with the tuning parameter $q_n$ is obtained as

$$\sqrt{n}\mathbf{V}_q^{-1/2}(\underline{\tilde{\theta}}_q - \underline{\theta}) \to N_p(\underline{0}, \mathbf{I}_p),$$

where $\mathbf{I}_p$ is the $(p \times p)$ identity matrix and

$$V_q(\underline{\theta}) = J(\underline{\theta}^*)^{-1} K(\underline{\theta}^*) J(\underline{\theta}^*)^{-1}$$

$$J(\underline{\theta}) = E\left(\nabla_{\underline{\theta}} U^*(\underline{\theta}_n^*; a_t, q_n)\right)$$

$$K(\underline{\theta}) = E\left(U^*(\underline{\theta}_n^*; a_t, q_n)' U^*(\underline{\theta}_n^*; a_t, q_n)\right).$$

Here $U^*(\underline{\theta}_n^*; a_t, q_n) = U(\underline{\theta}^*; a_t) f(a_t, \underline{\theta}^*)^{1-q_n}$ is the score vector and

$$\nabla_{\underline{\theta}} U^*(\underline{\theta}_n^*; a_t, q_n) = \nabla_{\underline{\theta}}\left(U(\underline{\theta}^*; a_t) f(a_t, \underline{\theta}^*)^{1-q_n}\right)$$

$$= \begin{cases} \nabla_{\underline{\theta}}\left(U(\underline{\theta}^*; a_t)\right) & , \quad q_n = 1 \\ f(a_t, \underline{\theta}^*)^{1-q_n}\left\{(1-q_n)U(\underline{\theta}^*; a_t)' U(\underline{\theta}^*; a_t) + \nabla_{\underline{\theta}} U(\underline{\theta}^*; a_t)\right\} & , \quad otherwise \end{cases}. \quad (22)$$

The first and second derivatives, which are required to obtain $K$ and $J$, are given in Appendix A and Appendix B.

A necessary and sufficient condition for asymptotic normality of CMLq estimator is $q_n \to 1$ when $n \to \infty$. It is important to note that when $q$ is fixed, CMLq estimator is equal to M-estimator, and so the asymptotic covariance matrix of CMLq estimator will be the same as the asymptotic covariance matrix of the M-estimator (Hampel, et al. 1986, Huber et al. 2009, Maronna et. al. 2006).

### 7. Numerical studies

In this section, to examine the performances of the CMLq method over the CML estimation method in the cases both with and without outliers in the data we carry on a simulation study and analyze a real-data. All computations are carried out in R-3.1.2 (R Development Core Team, 2017)

## 7.1 Sampling Design

We generate $y_t$ from model given in (3) included 5 and 10 covariates with same autoregressive structure $\underline{\phi} = (\phi_1, \phi_2)' = (0.8, -0.2)'$ and $x_{t,i}$ s are generated standard normal distribution. Here the values of $\underline{\phi}$ are taken to guarantee the stationarity assumption for the model of the error terms. Regression coefficients are considered as $\underline{\beta} = (\beta_1, \beta_2, \beta_3, \beta_4, \beta_5)' = (1, 3, 5, 2, 1)'$ and $\underline{\beta} = (\beta_1, \beta_2, ..., \beta_{10})' = (3, 3, ..., 3)'$. The simulation study is repeated 100 times with sample sizes 50.

The simulation study is reported under the following cases.

*Case I.* In the first part of our simulation study, data without contamination is considered. The results are presented in Table 1 and Table 4 for p=5 and p=10, respectively.

*Case II.* 10% of $y_i$ observations are replaced by the values generated from N(10,1), which are referred as outliers in y-direction. Table 2 and Table 5 show the simulation results for this scenario for p=5 and p=10, respectively.

*Case III.* In this case, both x and y variables will include outliers. In the y direction the outliers will be created using procedure described in Case II. In the x direction similar procedure will be used to create the outliers. Table 3 and Table 6 show the simulation results for this case for p=5 and p=10, respectively.

*Selection of the tuning parameter q*

The performance of the CMLq estimator depends on the tuning parameter that controls the weights. Therefore, it is important to select q appropriately.

In this paper we will use the robust AIC criterion proposed by Ronchetti (1985) to choose the appropriate tuning constant q. That is, we select the tuning parameter ($q$) which minimizes the following formula

$$RAIC = -\frac{1}{n}\sum_{t=p+1}^{n} L_q + tr(-M_2^{-1}M_1),$$

where $M_1$ is the first derivative and $M_2$ is the second derivative of $L_q$ function with respect to the parameters of interest. All of these derivatives are provided in Appendix A and Appendix B. The search is carried on using on the interval (0,1).

*Performance Measures*

To evaluate the performances, the bias, the root mean squared error (RMSE), the standard errors (SE) obtained from the asymptotic covariance matrix and the asymptotic confidence intervals (CIL − CIU) are calculated for each model parameter ($\gamma = \beta_i, \phi_l, \sigma$). The bias and RMSE are calculated using

$$Bias(\hat{\gamma}) = \bar{\gamma} - \gamma, \quad \bar{\gamma} = \frac{1}{100}\sum_{i=1}^{100}\hat{\gamma}_i,$$

$$RMSE(\hat{\gamma}) = \sqrt{\frac{1}{100}\sum_{i=1}^{100}(\hat{\gamma}_i - \gamma)^2}.$$

The asymptotic intervals for the CMLq estimates of the parameters of autoregressive error terms regression model are calculated by using the asymptotic covariance matrix given in Section 6.

Figures 1-6 are the boxplots of the model parameter estimates from 100 simulated datasets for the sample size 50 in all three cases.

## 7.2 Simulation Results

In Tables, root mean squared error (MSE), bias values, standard errors and the confidence intervals of the estimated parameters are given to compare the performance of the estimators. We also give the chosen $q$ values with respect to the minimum RAIC values.

The simulation results for Case I are summarized in Tables 1 for p=5 and in Tables 4 for p=10. Both Bias and RMSE values indicated that when there is no contamination in the data, the CMLq estimation is close to CML estimation method in every sample sizes. By examining the results came from different tuning parameter $q$, the best result for CMLq estimation is obtained from $q$ is close to 1 as expected. In Figure 1, Figure 4 we observe that the variability of the estimates obtained from the CML and CMLq estimation methods are very similar. Also, in those figures, boxplots show that both methods are accurately estimate the regression parameters.

Table 2 and Table 5 show the simulation results for Case II for p=5 and p=10, respectively. We observed from these results that the Bias and the RMSE values of the CMLq estimates are drastically better than those of the CML estimates. These results confirm that CML estimators are badly affected from the outliers in y direction.

The best results for CMLq estimator which are summarized in Table 5 are corresponding to the chosen q that makes the RAIC minimum. We observe that when the data contain outliers, the chosen q is not very close to one which makes the corresponding estimators robust against the outliers. Figure 2 and Figure 5 also show that if there are some outliers in the data the CML estimation method fails to correctly estimate the parameters. However, CMLq estimation method is not affected by the outliers in y direction and produces estimates that are very close to true parameter values. These results show that CMLq estimators are resistant to the vertical outliers in the data.

In the case of outliers in x-y direction, Table 3 and Table 6 show that CML estimation has a considerably higher RMSE than CMLq estimation for all the parameters except for the autoregressive model parameters. Figure 3 and Figure 6 can support the same results.

To sum up the simulation results confirm that the CMLq estimation method produces comparable results when there are no outliers in the data and it has definite superiority over the CML when there are outliers in the y and x and y direction. It should be noted that in case of outliers in x direction, the performance of the CMLq estimation method is not very premising. This is due to the fact that the CMLq estimation method is a M estimation method for fixed value of q and it is known that M-estimators are not robust against the outliers in x direction.

### 7.3 Real-data analysis

We consider a data set on the proportion of the number of ten million international phone calls from Belgium in the years 1950-1973. Rousseeuw and Leroy (1987) modeled this dataset with a robust regression method the least median of squares (LMS). Rousseeuw and Leroy (1987), point out that: "...it turned out that from 1964 to 1969 another recording system was used, giving the total number of minutes of these calls. The years 1963 and 1970 are also partially affected because the transitions did not happen exactly on New Year's Day...". This different measurement system in 1964-1969 causes a heavy contamination. As it can be seen in Figure 7, there are several outliers in the y-direction in 1964-1969.

Tuaç et al. (2017) have also used a linear regression model with AR(1) error terms with the assumption that the error terms have a t distribution as a heavy-tailed alternative to the normal distribution. For these reasons, this data set is modeled with autoregressive error terms regression model and the parameters of this model are estimated by using the CMLq method.

The number of phone calls made from Belgium is the dependent variable (y) and the explanatory variable (x) is the year. We observe from the autocorrelation function and the partial autocorrelation function graphs of the OLS residuals that the residuals show an autocorrelated structure with type AR(1). That is why we consider a regression model with AR(1) error term to model this data set. The summary of the results is reported in Table 7.

The linear regression model and the error structure are as follows

$$y_t = \beta_0 + \beta_1 x_t + e_t$$

$$e_t = \phi_1 + a_t$$

We also calculated RAIC for the real data example to compare the performance of the CMLq and CML estimation methods. The results are shown in Table 7. Figure 7 show the scatter plot of the data with the fitted regression lines obtained from CML and CMLq estimates. We observe from this figure that the CMLq estimate provides considerable better fit than the CML estimate.

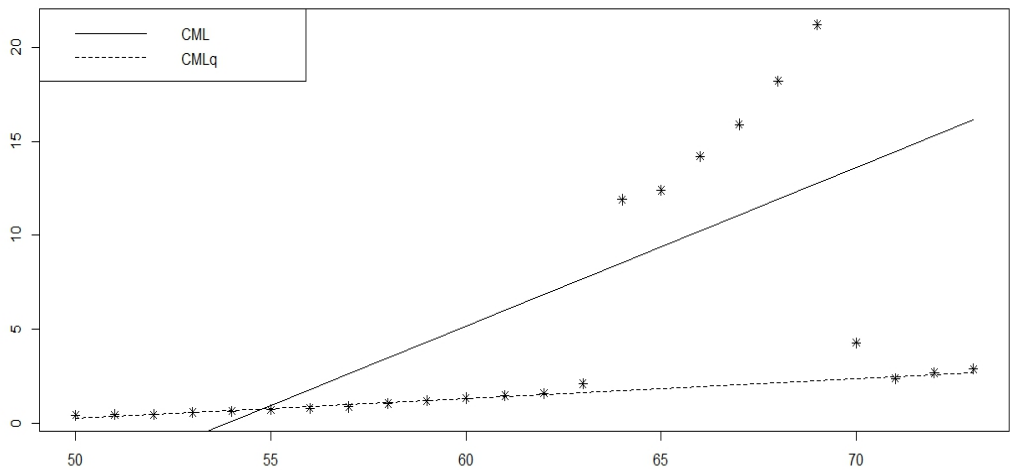

Figure 7. Number of international phone calls from Belgium with the CML and CMLq fit

Table 7. Summary of the estimated parameters for real data example

|  | CML | SE | CIL | CIU | CMLq | SE | CIL | CIU | q* |
|---|---|---|---|---|---|---|---|---|---|
| $\hat{\beta}_0$ | -45.36150 | 8.32749 | -48.96257 | -41.76042 | -5.0884 | 7.03634 | -8.13120 | -2.04571 | 0.917 |
| $\hat{\beta}_1$ | 0.84219 | 0,14078 | 0.78131 | 0.90307 | 0.10687 | 0.12175 | 0.05421 | 0.15952 |  |
| $\hat{\phi}$ | 0.81672 | 0.05827 | 0.79290 | 0.84053 | 0.70260 | 0.08389 | 0.66831 | 0.73688 |  |
| $\hat{\sigma}$ | 4.1424 | 0.08345 | 3.27576 | 5.99481 | 1.20222 | 0.04520 | 0.95068 | 1.73980 |  |
| RAIC | 61.3656 |  |  |  | 38.68 |  |  |  |  |

Also, according to results in Table 7, CMLq estimation method shows better performance in terms of the RAIC and the confidence intervals. We can observe that the length of the confidence intervals obtained from the CMLq are shorter than the length of the confidence intervals obtained from the CML estimators except the parameter $\phi$.

## 8 Discussion

In this paper, we have applied the CMLq estimation method to estimate the parameters of the regression model with autoregressive error term model when the data set is contaminated. The simulation results have revealed that CMLq estimation produces results that are close to the results obtained from the CML estimation method when there are no outliers in the data. However, the CMLq estimation outperforms when the dataset contains outliers in y and x and y direction. These simulation results have also displayed that the smaller values of q correspond to the CMLq estimators that are less sensitive to the outliers. To sum up, the CMLq estimation method can provide robust estimation method alternative to the CML and therefore this method can be used to estimate the parameters of the AR(p) error term regression model when outliers are present in the data.

**References**

Alpuim, T. and El-Shaarawi, A. 2008. On the efficiency of regression analysis with AR(p) errors. Journal of Applied Statistics, 35(7): 717-737.

Ansley, C. F. 1979. An algorithm for the exact likelihood of a mixed autoregressive-moving average process. Biometrika. 66 (1): 59-65.doi: 10.1093/biomet/66.1.59.

Beach, C.M. and Mackinnon, J. G. 1978. A Maximum Likelihood Procedure for Regression with Autocorrelated Errors. Econometrica, 46(1): 51- 58.

Cavalieri, J. 2002. O método de máxima Lq-verossimilhança em modelos com erros de medição (Doctoral thesis, Federal University of São Carlos, Department of Statistics). Retrieved from https://repositorio.ufscar.br/bitstream/handle/ufscar/4554/4180.pdf?sequence=1

Cochrane, D. and Orcutt, G. H. 1949. Application of least square to relationship containing autocorrelated error terms. Journal of American Statistical Association, 44: 32–61.

Dogru, F. Z., Bulut, Y. M. and Arslan, O. 2018. Doubly Reweighted Estimators for the Parameters of the Multivariate t Distribution. (Accepted)

Ferrari D. and Yang Y. 2007. Estimation of tail probability via the maximum Lq-likelihood method. Technical report 659, School of statistics, University of Minnesota

Ferrari D. and Paterlini S. 2009. The maximum Lq-likelihood method: an application to extreme quantile estimation in finance. Methodol Comput Appl Probab, 11(1): 3–19

Ferrari, D. and Yang, Y. 2010a. Maximum Lq-likelihood estimation. Annals of Statistics, 38(2): 753-783.

Ferrari, D., Paterlini, S. 2010b. Efficient and robust estimation for financial returns: an approach based on q-entropy. Available at SSRN: http://ssrn.com/abstract=1906819 or http://dx.doi.org/10.2139/ssrn.1906819.

Hampel, F. R., Ronchetti, E. M., Rousseeuw, P. J. and Stahel, W. A. 1986. Robust Statistics. The Approach Based on Influence Functions. New York: John Wiley and Sons. Chap. 4.2. pp. 230-231.

Havrda, J. and Charvát, F. 1967. Quantication method of classication processes: Concept of structural entropy. Kibernetika, 3: 30-35.

Huang, C., Lin, J. and Ren, Y.Y. 2013. Testing for the shape parameter of generalized extreme value distribution based on the Lq-likelihood ratio statistic. Metrika, 76: 641-671.

Huber, P. J. And Ronchetti, E. M. 2009. Robust Statistics. New Jersey: John Wiley and Sons.

Maronna, R. A., Martin, R. D. and Yohai, V. J. 2006. Robust Statistics: Theory and Methods. Chichester : John Wiley and Sons.

Ozdemir, S., Guney, Y. Tuac, Y. and Arslan, O. 2019. Maximum Lq-Likelihood Estimation for The Parameters of Marshall-Olkin Extended Burr XII Distribution. Commun. Fac. Sci. Univ. Ank. Ser. A1 Math. Stat, 68(1): 17-34.

Qin, Y. and Priebe, E.C. 2013. Maximum Lq-Likelihood Estimation via the Expectation Maximization Algorithm: A Robust Estimation of Mixture Models. Journal of the American Statistical Association, 108: 914-928.

Qin, Y. and Priebe, E.C. 2016. Robust Hypothesis Testing via Lq-Likelihood. Statistica Sinica: Preprint, doi:10.5705/ss.202015.0441.

R Core Team (2017). R: A language and environment for statistical computing. R Foundation for Statistical Computing, Vienna, Austria. URL http://www.R-project.org

Ronchetti, E. 1985. Robust model selection in regression. Statist. Probab. Lett, 3: 21–23.

Rousseeuw, P.J. and Leroy, A.M. 1987. Robust Regression and Outlier Detection. Wiley Series, USA.

Tsallis, C. 1988. Possible generalization of Boltzmann-Gibbs statistics. Journal of Statistical Physics, 52, 479-487.

Tuac, Y., Guney, Y., Senoglu, B. and Arslan, O. 2017. Robust Parameter Estimation of Regression Model with AR(p) Error Terms. Communications in Statistics - Simulation and Computation. https://doi.org/10.1080/03610918.2017.1343839 .

**APPENDIX A**

Let's $\underline{\theta} = (\beta_1, \beta_2, \ldots, \beta_M, \phi_1, \phi_2, \ldots, \phi_p, \sigma^2)$. The elements of $U(a_t; \underline{\theta})$ are

$$\frac{\partial \ln L}{\partial \beta_k} = \frac{1}{\sigma^2}\left(\Phi(B)y_t - \sum_{i=1}^{M} \beta_i \Phi(B)x_{t,i}\right)\Phi(B)x_{t,k}, \quad k = 1,2,\ldots,M$$

$$\frac{\partial \ln L}{\partial \phi_l} = \frac{1}{\sigma^2}\left(\Phi(B)y_t - \sum_{i=1}^{M} \beta_i \Phi(B)x_{t,i}\right)\left(y_{t-l} - \sum_{i=1}^{M} \beta_i x_{t-l,i}\right), \quad l = 1,2,\ldots,q,$$

$$\frac{\partial \ln L}{\partial \sigma^2} = -\frac{1}{2\sigma^2} + \frac{1}{2\sigma^4}\left(\Phi(B)y_t - \sum_{i=1}^{M} \beta_i \Phi(B)x_{t,i}\right)^2,$$

where $t = p+1, \ldots, N$.

The second partial derivatives of conditional log-likelihood function are

$$\frac{\partial^2 \ln L}{\partial \beta_j \partial \beta_k} = -\frac{1}{\sigma^2}\Phi(B)x_{t,j}\Phi(B)x_{t,k},$$

$$\frac{\partial^2 \ln L}{\partial \beta_j \partial \phi_i} = -\frac{1}{\sigma^2}\left[e_{t-i}\Phi(B)x_{t,j} + a_t x_{t-i,j}\right],$$

$$\frac{\partial^2 \ln L}{\partial \beta_j \partial \sigma^2} = -\frac{1}{\sigma^4}a_t \Phi(B)x_{t,j},$$

$$\frac{\partial^2 \ln L}{\partial \phi_i \partial \phi_r} = -\frac{1}{\sigma^2}e_{t-i}e_{t-r},$$

$$\frac{\partial^2 \ln L}{\partial \phi_i \partial \sigma^2} = -\frac{2}{\sigma^4}a_t e_{t-i},$$

$$\frac{\partial^2 \ln L}{\partial (\sigma^2)^2} = \frac{1}{2\sigma^4} - \frac{1}{\sigma^6}a_t^2,$$

where $t = p + 1, \ldots, N$, $j, k = 1, 2, \ldots, M$ and $i, r = 1, 2, \ldots, p$.

**APPENDIX B**

The first partial derivatives of Lq- likelihood function given in (), namely the elements of $U^*(a_t; \underline{\theta}, q)$ are

$$\frac{\partial Lq}{\partial \beta_k} = \frac{1}{\sigma^2} \omega_t \left( \Phi(B) y_t - \sum_{i=1}^{M} \beta_i \, \Phi(B) x_{t,i} \right) \Phi(B) x_{t,k}, \quad k = 1, 2, \ldots, M$$

$$\frac{\partial Lq}{\partial \phi_l} = \frac{1}{\sigma^2} \omega_t \left( \Phi(B) y_t - \sum_{i=1}^{M} \beta_i \, \Phi(B) x_{t,i} \right) \left( y_{t-l} - \sum_{i=1}^{M} \beta_i \, x_{t-l,i} \right), \quad l = 1, 2, \ldots, q,$$

$$\frac{\partial Lq}{\partial \sigma^2} = \omega_t \left[ -\frac{1}{2\sigma^2} + \frac{1}{2\sigma^4} \left( \Phi(B) y_t - \sum_{i=1}^{M} \beta_i \, \Phi(B) x_{t,i} \right)^2 \right],$$

where

$$\omega_t = \left( \frac{1}{\sqrt{2\pi\sigma^2}} \exp\left\{ -\frac{\left( \Phi(B) y_t - \sum_{i=1}^{M} \beta_i \, \Phi(B) x_{t,i} \right)^2}{2\sigma^2} \right\} \right)^{1-q}, t = p+1, \ldots N.$$

The second partial derivatives of conditional Lq-likelihood function can be obtained by substituting the derivatives given in Appendix A in the equation (22).

Table 1. Simulation results for Case I, p=5

| q | | Estimates | Bias | RMSE | SE | CIL | CIU |
|---|---|---|---|---|---|---|---|
| 0.99 | $\tilde{\beta}_1$ | 0.95568 | -0.04431 | 0.03271 | 0.032535 | 0.94666 | 0.96469 |
| | $\tilde{\beta}_2$ | 3.04248 | 0.04247 | 0.1836 | 0.186175 | 2.99087 | 3.09408 |
| | $\tilde{\beta}_3$ | 5.11161 | 0.11161 | 0.48159 | 0.48565 | 4.97699 | 5.24622 |
| | $\tilde{\beta}_4$ | 1.98414 | -0.01585 | 0.14032 | 0.137329 | 1.94607 | 2.02220 |
| | $\tilde{\beta}_5$ | 0.91404 | -0.08596 | 0.33923 | 0.340653 | 0.81961 | 1.00846 |
| | $\tilde{\phi}_1$ | 0.85362 | 0.05362 | 0.11338 | 0.19614 | 0.79925 | 0.90798 |
| | $\tilde{\phi}_2$ | -0.23565 | -0.03565 | 0.14057 | 0.18256 | -0.28625 | -0.18505 |
| | $\tilde{\sigma}$ | 0.95734 | -0.04266 | 0.11086 | 0.16592 | 0.82682 | 1.25555 |
| 1.00 | $\hat{\beta}_1$ | 0.94579 | -0.0542 | 0.04062 | 0.04011 | 0.934672 | 0.95690 |
| | $\hat{\beta}_2$ | 3.08820 | 0.0882 | 0.33399 | 0.33698 | 2.994794 | 3.18160 |
| | $\hat{\beta}_3$ | 5.1055 | 0.10549 | 0.43553 | 0.4364 | 4.984536 | 5.22646 |
| | $\hat{\beta}_4$ | 1.97433 | -0.02566 | 0.14209 | 0.1403 | 1.935441 | 2.01321 |
| | $\hat{\beta}_5$ | 0.95916 | -0.04083 | 0.33732 | 0.33845 | 0.865346 | 1.05297 |
| | $\hat{\phi}_1$ | 0.83571 | 0.0357 | 0.17236 | 0.17432 | 0.787391 | 0.88402 |
| | $\hat{\phi}_2$ | -0.22147 | -0.02146 | 0.13431 | 0.17758 | -0.27069 | -0.17225 |
| | $\hat{\sigma}$ | 0.96519 | -0.0348 | 0.10839 | 0.12213 | 0.82009 | 1.24533 |

Table 2. Simulation results for Case II, p=5

| q | | Estimates | Bias | RMSE | SE | CIL | CIU |
|---|---|---|---|---|---|---|---|
| 0.60 | $\tilde{\beta}_1$ | 1.04486 | 0.04486 | 0.29874 | 0.29647 | 0.96268 | 1.12703 |
| | $\tilde{\beta}_2$ | 2.89744 | -0.10255 | 0.47066 | 0.47394 | 2.76607 | 3.02880 |
| | $\tilde{\beta}_3$ | 5.00725 | 0.00725 | 0.04656 | 0.04251 | 4.99546 | 5.01903 |
| | $\tilde{\beta}_4$ | 2.08757 | 0.08757 | 0.37546 | 0.38353 | 1.98126 | 2.19387 |
| | $\tilde{\beta}_5$ | 1.0461 | 0.04610 | 0.31289 | 0.30953 | 0.96030 | 1.13189 |
| | $\tilde{\phi}_1$ | 0.76276 | -0.03723 | 0.18873 | 0.18428 | 0.71168 | 0.81384 |
| | $\tilde{\phi}_2$ | -0.04947 | 0.15052 | 0.64675 | 0.75916 | -0.25990 | 0.16095 |
| | $\tilde{\sigma}$ | 1.01113 | 0.01113 | 0.51718 | 0.52521 | 0.86618 | 1.31531 |
| 1.00 | $\hat{\beta}_1$ | 1.43832 | 0.43832 | 1.98659 | 1.93437 | 0.90214 | 1.97450 |
| | $\hat{\beta}_2$ | 4.06013 | 1.06014 | 3.845018 | 4.45256 | 2.82594 | 5.29431 |
| | $\hat{\beta}_3$ | 6.75896 | 1.75897 | 4.256039 | 6.91579 | 4.84200 | 8.67591 |
| | $\hat{\beta}_4$ | 2.48166 | 0.48166 | 1.209248 | 2.58094 | 1.76626 | 3.19706 |
| | $\hat{\beta}_5$ | 1.70155 | 0.70155 | 1.452112 | 2.67861 | 0.95907 | 2.44402 |
| | $\hat{\phi}_1$ | -0.11389 | -0.91389 | 1.156451 | 3.48156 | -1.07893 | 0.85114 |
| | $\hat{\phi}_2$ | 0.05842 | 0.25843 | 0.981569 | 1.05701 | -0.23457 | 0.35140 |
| | $\hat{\sigma}$ | 10.23524 | 9.23524 | 10.06469 | 9.28857 | 0.76790 | 13.31430 |

Table 3. Simulation results for Case III, p=5

| q | | Estimates | Bias | RMSE | SE | CIL | CIU |
|---|---|---|---|---|---|---|---|
| 0.44 | $\tilde{\beta}_1$ | 0.97228 | -0.02771 | 0.36524 | 0.35617 | 0.87355 | 1.07100 |
| | $\tilde{\beta}_2$ | 3.10851 | 0.10851 | 0.39046 | 0.48200 | 2.97490 | 3.24211 |
| | $\tilde{\beta}_3$ | 4.93869 | -0.0613 | 0.50891 | 0.49178 | 4.80237 | 5.07500 |
| | $\tilde{\beta}_4$ | 2.08666 | 0.08666 | 0.31819 | 0.40511 | 1.97436 | 2.19895 |
| | $\tilde{\beta}_5$ | 0.89745 | -0.10254 | 0.52575 | 0.48494 | 0.76303 | 1.03186 |
| | $\tilde{\phi}_1$ | 0.79471 | -0.00528 | 0.09552 | 0.09315 | 0.76889 | 0.82053 |
| | $\tilde{\phi}_2$ | -0.04648 | 0.15351 | 0.64144 | 0.64641 | -0.22566 | 0.13269 |
| | $\tilde{\sigma}$ | 1.11441 | 0.11441 | 0.6555 | 0.6941 | 0.95465 | 1.44966 |
| 1.00 | $\hat{\beta}_1$ | 0.55661 | -0.44338 | 1.77338 | 1.93625 | 0.01990 | 1.09331 |
| | $\hat{\beta}_2$ | 1.8251 | -1.1749 | 2.00342 | 4.26081 | 0.64406 | 3.00613 |
| | $\hat{\beta}_3$ | 2.81999 | -2.18001 | 5.43887 | 8.40224 | 0.49100 | 5.14897 |
| | $\hat{\beta}_4$ | 1.26058 | -0.73941 | 1.69692 | 2.9974 | 0.42974 | 2.09141 |
| | $\hat{\beta}_5$ | 0.45004 | -0.54996 | 1.81343 | 2.94925 | -0.36745 | 1.26753 |
| | $\hat{\phi}_1$ | 0.06648 | -0.73352 | 1.77549 | 2.81078 | -0.71263 | 0.84558 |
| | $\hat{\phi}_2$ | -0.18028 | 0.01972 | 0.6972 | 0.65814 | -0.36271 | 0.00214 |
| | $\hat{\sigma}$ | 9.52381 | 8.52382 | 8.99814 | 1.75666 | 0.62334 | 12.38886 |

Table 4. Simulation results for Case I, p=10

| q | | Estimates | Bias | RMSE | SE | CIL | CIU |
|---|---|---|---|---|---|---|---|
| 0.99 | $\tilde{\beta}_1$ | 2.99046 | -0.00954 | 0.09088 | 0.08861 | 2.96589 | 3.01502 |
| | $\tilde{\beta}_2$ | 2.94312 | -0.05688 | 0.19870 | 0.30373 | 2.85893 | 3.02731 |
| | $\tilde{\beta}_3$ | 2.96989 | -0.03011 | 0.18171 | 0.19418 | 2.91606 | 3.02371 |
| | $\tilde{\beta}_4$ | 2.99039 | -0.00961 | 0.08300 | 0.07931 | 2.96840 | 3.01237 |
| | $\tilde{\beta}_5$ | 2.98233 | -0.01767 | 0.12804 | 0.12564 | 2.94750 | 3.01715 |
| | $\tilde{\beta}_6$ | 2.98382 | -0.01618 | 0.10177 | 0.10730 | 2.95407 | 3.01356 |
| | $\tilde{\beta}_7$ | 2.99666 | -0.00334 | 0.07957 | 0.07999 | 2.97448 | 3.01883 |
| | $\tilde{\beta}_8$ | 2.99154 | -0.00846 | 0.08836 | 0.08791 | 2.96717 | 3.01590 |
| | $\tilde{\beta}_9$ | 2.98169 | -0.01831 | 0.12315 | 0.12893 | 2.94595 | 3.01742 |
| | $\tilde{\beta}_{10}$ | 3.06193 | 0.06193 | 0.31800 | 0.31167 | 2.97553 | 3.14832 |
| | $\tilde{\phi}_1$ | 0.78136 | -0.01864 | 0.12322 | 0.12308 | 0.74724 | 0.81547 |
| | $\tilde{\phi}_2$ | -0.20134 | -0.00134 | 0.04165 | 0.0423 | -0.21306 | -0.18962 |
| | $\tilde{\sigma}$ | 1.00999 | 0.00999 | 0.00876 | 0.00881 | 0.90836 | 1.41563 |
| 1.00 | $\hat{\beta}_1$ | 2.99532 | -0.00468 | 0.08804 | 0.08732 | 2.97111 | 3.01952 |
| | $\hat{\beta}_2$ | 2.94274 | -0.05726 | 0.25417 | 0.25364 | 2.87243 | 3.01304 |
| | $\hat{\beta}_3$ | 2.97302 | -0.02698 | 0.09807 | 0.10266 | 2.94456 | 3.00147 |
| | $\hat{\beta}_4$ | 3.00258 | 0.00258 | 0.07348 | 0.07215 | 2.98258 | 3.02257 |
| | $\hat{\beta}_5$ | 2.98221 | -0.01779 | 0.13158 | 0.12087 | 2.94870 | 3.01571 |
| | $\hat{\beta}_6$ | 2.98584 | -0.01416 | 0.09284 | 0.09489 | 2.95953 | 3.01214 |
| | $\hat{\beta}_7$ | 2.99757 | -0.00243 | 0.07611 | 0.07328 | 2.97725 | 3.01788 |
| | $\hat{\beta}_8$ | 2.99822 | -0.00178 | 0.07573 | 0.07837 | 2.97649 | 3.01994 |
| | $\hat{\beta}_9$ | 2.98818 | -0.01182 | 0.11985 | 0.12292 | 2.95410 | 3.02225 |
| | $\hat{\beta}_{10}$ | 3.05098 | 0.05098 | 0.2166 | 0.23053 | 2.98708 | 3.11488 |
| | $\hat{\phi}_1$ | 0.77994 | -0.02006 | 0.12453 | 0.12090 | 0.74642 | 0.81345 |
| | $\hat{\phi}_2$ | -0.21621 | -0.01621 | 0.04228 | 0.04599 | -0.22896 | -0.20346 |
| | $\hat{\sigma}$ | 1.01461 | 0.01461 | 0.00913 | 0.08417 | 0.91251 | 1.42210 |

Table 5. Simulation results for Case II, p=10

| q | | Estimates | Bias | RMSE | SE | CIL | CIU |
|---|---|---|---|---|---|---|---|
| 0.63 | $\tilde{\beta}_1$ | 3.10437 | 0.10437 | 0.46401 | 0.46874 | 2.97444 | 3.23430 |
| | $\tilde{\beta}_2$ | 2.99613 | -0.00386 | 0.11980 | 0.12620 | 2.96115 | 3.03111 |
| | $\tilde{\beta}_3$ | 3.17695 | 0.17695 | 0.93522 | 0.98024 | 2.90524 | 3.44866 |
| | $\tilde{\beta}_4$ | 2.99576 | -0.00423 | 0.06135 | 0.06784 | 2.97696 | 3.01456 |
| | $\tilde{\beta}_5$ | 2.89541 | -0.10458 | 0.45475 | 0.44097 | 2.77318 | 3.01764 |
| | $\tilde{\beta}_6$ | 2.89955 | -0.10044 | 0.62792 | 0.61124 | 2.73012 | 3.06898 |
| | $\tilde{\beta}_7$ | 3.21677 | 0.21677 | 0.78716 | 0.92172 | 2.96128 | 3.47226 |
| | $\tilde{\beta}_8$ | 2.95599 | -0.04400 | 0.31494 | 0.34858 | 2.85937 | 3.05261 |
| | $\tilde{\beta}_9$ | 2.80923 | -0.19076 | 0.88541 | 0.98674 | 2.53572 | 3.08274 |
| | $\tilde{\beta}_{10}$ | 3.35013 | 0.45013 | 0.97399 | 1.26457 | 2.99961 | 3.70065 |
| | $\tilde{\phi}_1$ | 0.76617 | -0.03382 | 0.23361 | 0.54672 | 0.61463 | 0.91771 |
| | $\tilde{\phi}_2$ | 0.06148 | 0.26148 | 0.80275 | 0.97569 | -0.20897 | 0.33193 |
| | $\tilde{\sigma}$ | 1.19711 | 0.19711 | 0.82447 | 0.73086 | 0.97665 | 2.67790 |
| 1.00 | $\hat{\beta}_1$ | 3.32414 | 0.32415 | 1.34633 | 1.21082 | 2.98852 | 3.65976 |
| | $\hat{\beta}_2$ | 3.41139 | 0.41140 | 1.92714 | 2.01287 | 2.85345 | 3.96933 |
| | $\hat{\beta}_3$ | 4.17681 | 1.17682 | 1.47198 | 4.64121 | 2.89033 | 5.46329 |
| | $\hat{\beta}_4$ | 4.31203 | 1.31203 | 3.25816 | 5.26659 | 2.85221 | 5.77185 |
| | $\hat{\beta}_5$ | 3.94000 | 0.94001 | 1.23295 | 3.71485 | 2.91030 | 4.96970 |
| | $\hat{\beta}_6$ | 4.44877 | 1.44877 | 3.40305 | 5.80092 | 2.84084 | 6.05670 |
| | $\hat{\beta}_7$ | 3.69657 | 0.69657 | 2.94512 | 3.18861 | 2.81273 | 4.58041 |
| | $\hat{\beta}_8$ | 3.97170 | 0.97171 | 2.89445 | 3.62755 | 2.96619 | 4.97721 |
| | $\hat{\beta}_9$ | 3.88432 | 0.88433 | 2.92778 | 3.83186 | 2.82218 | 4.94646 |
| | $\hat{\beta}_{10}$ | 4.04412 | 1.04413 | 2.06263 | 3.97276 | 2.94293 | 5.14531 |
| | $\hat{\phi}_1$ | 0.00740 | -0.79260 | 1.02508 | 2.99382 | -0.82244 | 0.83724 |
| | $\hat{\phi}_2$ | -0.38703 | -0.18703 | 0.81849 | 0.79606 | -0.60769 | -0.16637 |
| | $\hat{\sigma}$ | 13.22522 | 12.22523 | 13.27970 | 12.03498 | 0.28287 | 18.53680 |

Table 6. Simulation results for Case III, p=10

| q | | Estimates | Bias | RMSE | SE | CIL | CIU |
|---|---|---|---|---|---|---|---|
| 0.41 | $\tilde{\beta}_1$ | 2.87524 | -0.12475 | 0.34051 | 0.48154 | 2.74176 | 3.00872 |
| | $\tilde{\beta}_2$ | 2.65257 | -0.34742 | 0.91942 | 1.30643 | 2.29045 | 3.01469 |
| | $\tilde{\beta}_3$ | 2.81819 | -0.18181 | 0.6982 | 0.81535 | 2.59219 | 3.04419 |
| | $\tilde{\beta}_4$ | 2.77023 | -0.22977 | 0.642591 | 0.84573 | 2.53581 | 3.00465 |
| | $\tilde{\beta}_5$ | 2.80791 | -0.19209 | 0.93393 | 0.93361 | 2.54913 | 3.06669 |
| | $\tilde{\beta}_6$ | 2.76750 | -0.23249 | 0.850986 | 0.87468 | 2.52505 | 3.00995 |
| | $\tilde{\beta}_7$ | 2.98588 | -0.01412 | 0.23973 | 0.21890 | 2.92520 | 3.04656 |
| | $\tilde{\beta}_8$ | 2.82776 | -0.17223 | 0.647294 | 0.68577 | 2.63767 | 3.01785 |
| | $\tilde{\beta}_9$ | 2.74666 | -0.25333 | 0.96672 | 0.97839 | 2.47546 | 3.01786 |
| | $\tilde{\beta}_{10}$ | 2.85908 | -0.14092 | 0.650834 | 0.70445 | 2.66382 | 3.05434 |
| | $\tilde{\phi}_1$ | 0.67642 | -0.12358 | 0.50757 | 0.48179 | 0.54287 | 0.80997 |
| | $\tilde{\phi}_2$ | 0.03722 | 0.23722 | 0.961066 | 0.96167 | -0.22934 | 0.30378 |
| | $\tilde{\sigma}$ | 2.35923 | 1.35923 | 1.62963 | 1.96340 | 0.52184 | 3.30676 |
| 1.00 | $\hat{\beta}_1$ | 2.31931 | -0.68069 | 1.37979 | 3.07655 | 1.46653 | 3.17209 |
| | $\hat{\beta}_2$ | 2.23790 | -0.76210 | 1.17922 | 2.85772 | 1.44578 | 3.03002 |
| | $\hat{\beta}_3$ | 2.28866 | -0.71134 | 2.34043 | 3.09404 | 1.43104 | 3.14628 |
| | $\hat{\beta}_4$ | 2.24767 | -0.75233 | 2.24115 | 2.82917 | 1.46346 | 3.03188 |
| | $\hat{\beta}_5$ | 2.34417 | -0.65583 | 1.15104 | 2.75168 | 1.58144 | 3.10690 |
| | $\hat{\beta}_6$ | 2.29419 | -0.70581 | 1.29847 | 2.62589 | 1.56633 | 3.02205 |
| | $\hat{\beta}_7$ | 2.43387 | -0.56613 | 1.14423 | 2.29822 | 1.79684 | 3.07090 |
| | $\hat{\beta}_8$ | 2.05707 | -0.94293 | 1.39418 | 3.63473 | 1.04957 | 3.06457 |
| | $\hat{\beta}_9$ | 2.16312 | -0.83688 | 1.39916 | 3.18660 | 1.27984 | 3.04640 |
| | $\hat{\beta}_{10}$ | 2.29770 | -0.70230 | 1.10200 | 2.85404 | 1.50660 | 3.08880 |
| | $\hat{\phi}_1$ | 0.09797 | -0.70203 | 0.91911 | 2.78527 | -0.67407 | 0.87001 |
| | $\hat{\phi}_2$ | -0.15374 | 0.04626 | 0.58027 | 0.61885 | -0.32528 | 0.01780 |
| | $\hat{\sigma}$ | 7.09458 | 6.09458 | 6.62821 | 6.35406 | 0.38068 | 9.94393 |

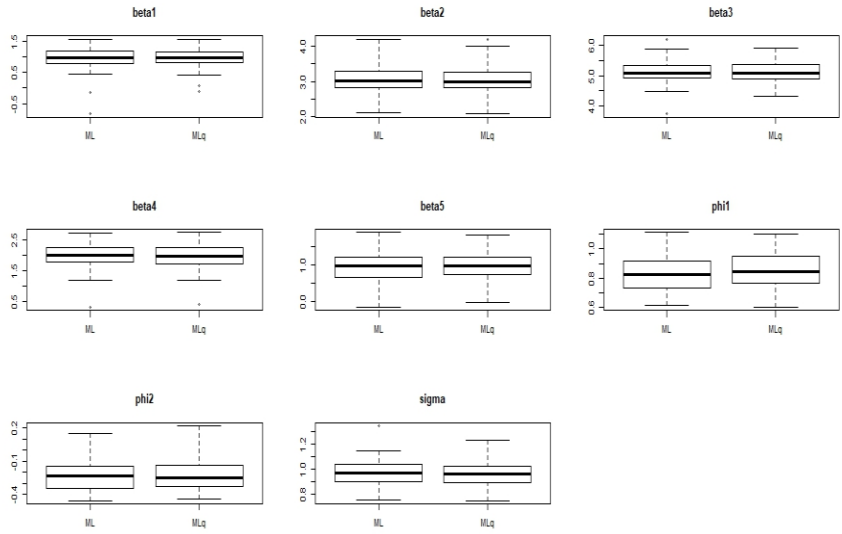

Figure 1. Boxplots of the estimates of the parameters for Case I, p=5. True parameter values are $\underline{\beta} = (1, 3, 5, 2, 1)'$, $\underline{\phi} = (0.8, -0.2)'$, $\sigma = 1$.

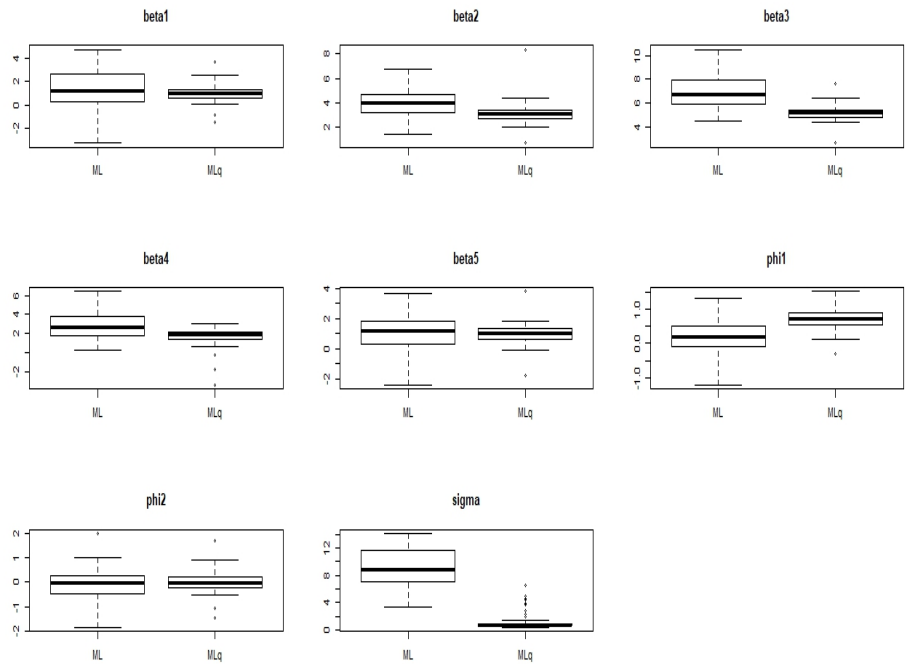

Figure 2. Boxplots of the estimates of the parameters for Case II, p=5. True parameter values are $\underline{\beta} = (1, 3, 5, 2, 1)'$, $\underline{\phi} = (0.8, -0.2)'$, $\sigma = 1$

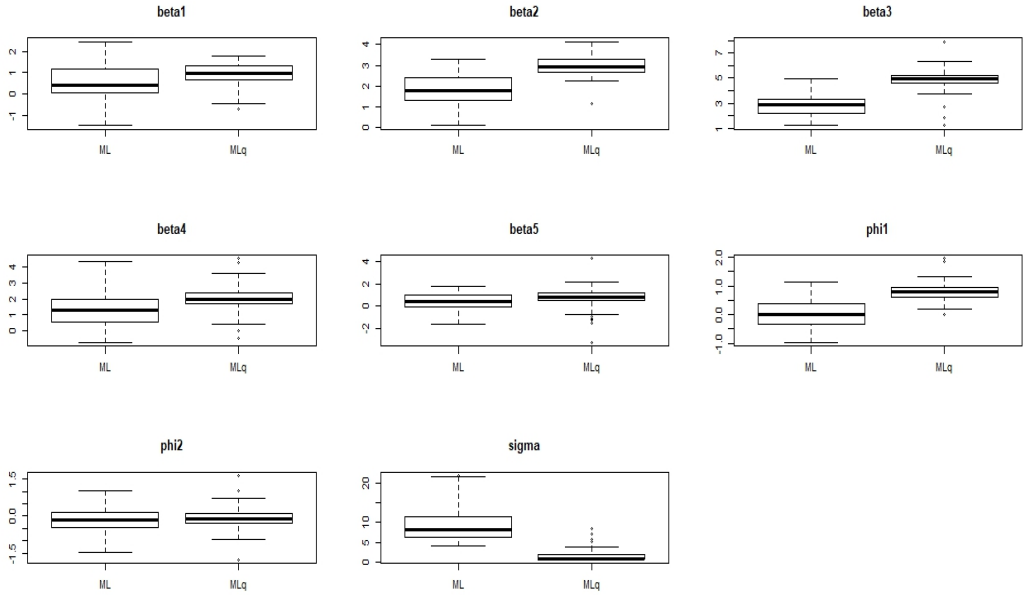

Figure 3. Boxplots of the estimates of the parameters for Case III, p=5. True parameter values are $\underline{\beta} = (1, 3, 5, 2, 1)'$, $\underline{\phi} = (0.8, -0.2)'$, $\sigma = 1$.

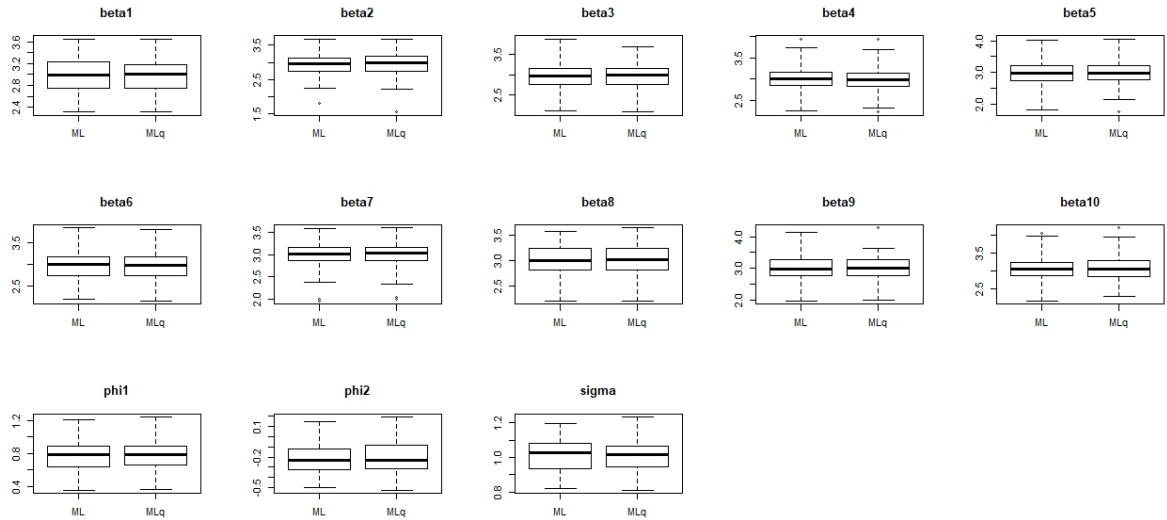

Figure 4. Boxplots of the estimates of the parameters for Case I, p=10. True parameter values are $\underline{\beta} = (\beta_1, \beta_2, \ldots, \beta_{10})' = (3, 3, \ldots, 3)'$, $\underline{\phi} = (0.8, -0.2)'$, $\sigma = 1$.

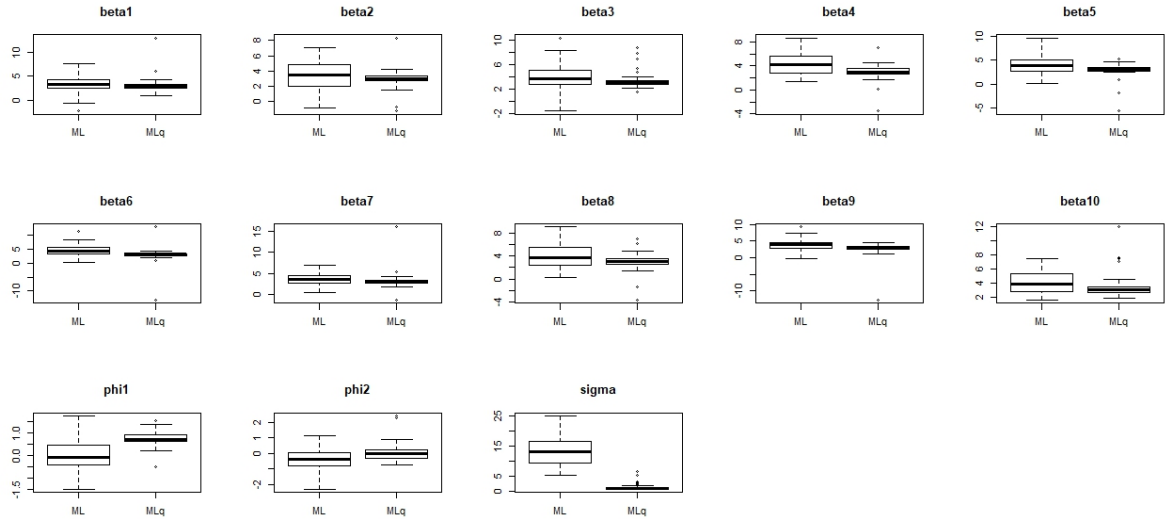

Figure 5. Boxplots of the estimates of the parameters for Case II, p=10. True parameter values are $\underline{\beta} = (\beta_1, \beta_2, \ldots, \beta_{10})' = (3, 3, \ldots, 3)'$, $\underline{\phi} = (0.8, -0.2)'$, $\sigma = 1$.

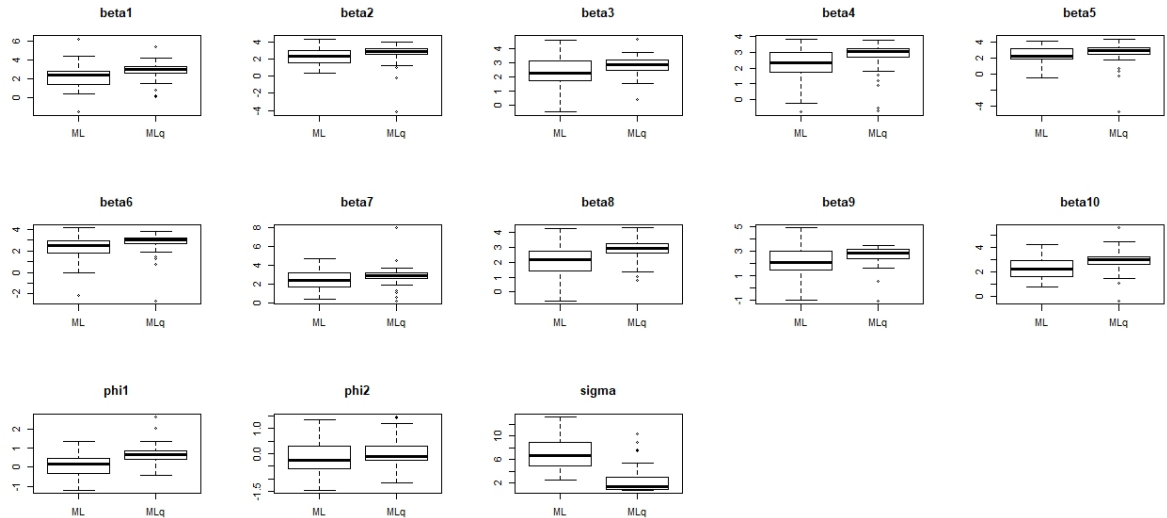

Figure 6. Boxplots of the estimates of the parameters for Case III, p=10. True parameter values are $\underline{\beta} = (\beta_1, \beta_2, \ldots, \beta_{10})' = (3, 3, \ldots, 3)'$, $\underline{\phi} = (0.8, -0.2)'$, $\sigma = 1$.

35